\documentclass[10pt]{amsart}
\usepackage{graphicx}
\usepackage{amsmath}
\usepackage{amssymb}
\usepackage[margin=2cm]{geometry}
\usepackage{hyperref}
\usepackage{algorithm}
\usepackage{algpseudocode}
\algnewcommand\algorithmicinput{\textbf{Input:}}
\algnewcommand\algorithmicoutput{\textbf{Output:}}
\algnewcommand\Input{\item[\algorithmicinput]}
\algnewcommand\Output{\item[\algorithmicoutput]}\usepackage{caption} % per \captionof
\usepackage{booktabs}
\usepackage{longtable}
\usepackage{tabularx}
\usepackage[table]{xcolor}
\usepackage{placeins}
\usepackage{amsthm}
\usepackage[backend=biber]{biblatex}
\newtheorem{remark}{Remark}

\title[A reliability-aware randomized simheuristic for the stochastic team orienteering problem]{A reliability-aware randomized simheuristic for the stochastic team orienteering problem}

\author{Michele Circelli}

\begin{document}

\begin{abstract}
We study a stochastic variant of the Team Orienteering Problem with lognormal travel
times and an all-or-nothing reward policy, under which the reward of a route is lost if
its travel time exceeds the available budget. We propose a reliability-aware simheuristic
that combines a savings-based constructive heuristic with three specific design elements:
a Top-$L_{\mathrm{top}}$ randomization mechanism, stochastic screening of the savings
parameter, and an explicit reliability threshold for solution selection. Computational
experiments on the Chao et al.\ benchmark show that the method is competitive with the
VNS-based simheuristic of Panadero et al.\ (2020) on a non-trivial subset of instances
using a significantly simpler architecture, with the largest gains on two-vehicle
sub-families with long routes where the reliability-aware selection compensates for the
absence of VNS-style exploration. On larger multi-vehicle instances the simpler
architecture is outperformed, and this trade-off is discussed explicitly.
\end{abstract}

\maketitle

\section*{Introduction}

The Team Orienteering Problem (TOP) is a well-known variant of vehicle routing in which
a fleet of vehicles must select and visit a subset of nodes within a limited time budget
to maximize total collected reward. Due to these characteristics, the TOP has been
successfully applied to service routing, surveillance, inspection planning, and related
contexts~\cite{Chao1,gunawan2016orienteering,caceres2015rich}.

In many practical settings travel times are uncertain, motivating the study of the
Stochastic Team Orienteering Problem (STOP), where travel times are modeled as random
variables. In this setting, solutions must balance expected reward against the risk of
violating route-duration constraints. Deterministic solutions that appear attractive
under nominal conditions may become fragile under uncertainty, losing both expected
reward and reliability~\cite{campbell2011orienteering,juan2015simheuristics}. Despite
its practical relevance, the stochastic TOP has received much less attention than its
deterministic counterpart, and fewer works address the multi-vehicle case with an
all-or-nothing reward structure.

This paper proposes a reliability-aware simheuristic for the STOP with lognormal travel
times and all-or-nothing route rewards. The method builds on the savings-based
construction pipeline of Panadero et al.~\cite{Juan}, which is adopted as the
deterministic backbone without modification. The three algorithmic contributions of this
paper are:
\begin{enumerate}
    \item \textbf{Top-$L_{\mathrm{top}}$ randomization}: moves are drawn uniformly among
    the top $L_{\mathrm{top}}$ candidates, providing a simpler and more transparent
    alternative to geometric biased randomization with a single interpretable parameter.
    \item \textbf{Stochastic screening of $\alpha$}: candidate values of the savings
    parameter are evaluated directly under stochastic travel times rather than being
    fixed on a deterministic criterion, ensuring consistency with the stochastic
    objective throughout.
    \item \textbf{Reliability-aware selection}: solutions are filtered by a minimum
    reliability threshold $\beta$ before being compared by expected reward, formalizing
    the trade-off between robustness and performance under the all-or-nothing policy.
\end{enumerate}
These components are simpler and more modular than the VNS-based machinery of prior
work, and are designed to be transferable to other stochastic combinatorial optimization
problems.

The remainder of the paper is organized as follows. The literature review summarizes the
most relevant related contributions. Section~\ref{mathsubsec} presents the mathematical
model. Section~\ref{algorithmsubsec} describes the algorithmic framework.
Section~\ref{sec:results} reports the computational experiments.
Section~\ref{sec:conclusions} concludes and outlines future research.

\subsection*{Literature Review}\label{sec:literature}

The Orienteering Problem (OP) is a classical routing problem in which a vehicle must
select a subset of customers to visit within a limited budget to maximize collected
reward. Since its introduction by Golden et al.~\cite{golden1987orienteering}, the OP
and its variants have attracted substantial attention; surveys can be found in
Vansteenwegen et al.~\cite{vansteenwegen2011orienteering} and Gunawan
et al.~\cite{gunawan2016orienteering}, while C\'aceres-Cruz et al.~\cite{caceres2015rich}
discuss it within the wider class of rich vehicle routing problems. Chao
et al.~\cite{Chao2} proposed an early effective heuristic for the single-vehicle case.

The Team Orienteering Problem~\cite{Chao1} extends the OP to multiple vehicles. Because
the TOP is NP-hard, exact methods are limited to medium-sized instances --- Butt and
Ryan~\cite{butt1999} proposed a column-generation approach for moderate-size problems ---
and the literature has mainly focused on metaheuristics. These include tabu search and
VNS~\cite{archetti2007,tang2005b}, particle swarm optimization~\cite{dang2013},
multi-start simulated annealing~\cite{lin2013}, and evolutionary
approaches~\cite{ferreira2014,ke2016}.

Stochastic variants of the TOP have received much less attention. Early work on
stochastic orienteering mainly focused on the single-vehicle setting, with contributions
addressing uncertainty in rewards~\cite{ilhan2008}, stochastic travel and service
times~\cite{campbell2011orienteering,papapanagiotou2014,evers2014,verbeeck2016}, and
chance-constrained or recourse-based
formulations~\cite{tang2005a,lau2012,varakantham2013,zhang2014}.

Simheuristics --- which combine metaheuristic search with Monte Carlo simulation ---
have emerged as a flexible approach for stochastic combinatorial optimization. Juan
et al.~\cite{juan2015simheuristics} provide a comprehensive review; applications include
stochastic routing~\cite{gruler2018,gruler2017,gonzalezmartin2018}. Within this
paradigm, biased randomization~\cite{grasas2017biased} provides a way to diversify
constructive heuristics while preserving computational efficiency.

The application of simheuristics to the STOP was initiated by Panadero
et al.~\cite{panadero2017stop} and developed into a full contribution in~\cite{Juan},
which proposed a savings-based heuristic embedded in a VNS framework with Monte Carlo
simulation under an all-or-nothing reward policy. Subsequent work extended this line to
dynamic rewards~\cite{bayliss2020}, probabilistic
delays~\cite{herrera2022}, position-dependent rewards~\cite{panadero2022position},
dynamic and learning-based settings~\cite{peyman2024,uguina2024}, and real-time
methods~\cite{panadero2021applsci,reyesrubiano2020dynamic}. Panadero
et al.~\cite{panadero2024itor} presented a comparative study benchmarking a BRVNS simheuristic
against the SAA method and a hybrid model, showing that the simheuristic outperforms
both in solution quality and computational time. From a different perspective, Yu
et al.~\cite{yu2022robust} addressed a robust TOP with decreasing profits using
distributionally robust optimization.

The present paper addresses the same STOP setting as~\cite{Juan,panadero2024itor} --- lognormal
travel times, all-or-nothing policy, Chao benchmark --- and contributes a specific
algorithmic design distinct from the VNS-based and biased-randomized approaches of prior
work: Top-$L_{\mathrm{top}}$ randomization, stochastic screening of $\alpha$, and
reliability-aware selection via $\beta$, together with an empirical assessment of the
resulting simplicity--performance trade-off.

\subsection*{Problem statement}\label{sec:problem}

We consider a routing problem in which a fleet of $m$ vehicles must visit a subset of
geographically dispersed nodes within a limited operational time budget $T_{\max} > 0$.
Each node $i$ has a nonnegative reward $u_i \geq 0$ and can be visited at most once.
Every vehicle departs from a start depot and terminates at an end depot.

The objective is to determine at most $m$ start-to-end routes maximizing total collected
reward. In the deterministic setting, a route is feasible if its total travel time does
not exceed $T_{\max}$. In the stochastic setting, travel times are random variables. Under
the all-or-nothing policy, if the realized travel time of a vehicle exceeds $T_{\max}$,
all rewards associated with that route are lost. Uncertainty therefore affects both the
feasibility of routes and the expected value of a solution.

We address the stochastic version of the problem, in which travel times are modeled as log-normal random variables and the objective is to maximize the expected collected reward under an all-or-nothing reward policy. A deterministic formulation — with Euclidean travel times — is also introduced, but solely as the structural backbone of the proposed simheuristic; no deterministic benchmark is reported and no contribution is claimed on that side. Solution reliability, defined as the probability that all routes complete within the time budget, plays a central role in the evaluation and selection of candidate solutions.

\section{Mathematical model}\label{mathsubsec}

In this section, we introduce a MILP formulation of the deterministic problem and its
stochastic extension.

We model the problem on a complete directed graph. Let
\[
N=\{0,1,\dots,n,n+1\}
\]
be the set of nodes, where $0$ denotes the start depot, nodes $i=1,\dots,n$ represent
the customer nodes, and $n+1$ denotes the end depot, with $n\in\mathbb{N}\setminus\{0\}$.

Each node $i\in N$ is associated with planar coordinates $(x_i,y_i)\in\mathbb{R}^2$.
Each customer node $i\in \{1,\dots,n\}$ has a nonnegative reward $u_i\ge 0$, while we
set $u_0:=u_{n+1}:=0$.

The arc set is
\[
A := \{(i,j)\in N\times N : i\neq j\},
\]
and we denote by
\[
a:=|A|=(n+2)(n+1).
\]

Let $M=\{1,\dots,m\}$ be the set of available vehicles, with $m\in\mathbb{N}\setminus\{0\}$.
Each vehicle has a maximum travel-time budget $T_{\max}>0$.

\subsection{Decision variables}

For each vehicle $k\in M$ and each arc $(i,j)\in A$, we introduce binary routing
variables $\rho_{ij}^k\in\{0,1\}$:
\[
\rho_{ij}^k=
\begin{cases}
1, &\text{if vehicle $k$ travels directly from node $i$ to node $j$},\\
0, &\text{otherwise.}
\end{cases}
\]
We also introduce binary activation variables $\delta_k\in\{0,1\}$ for each vehicle
$k\in M$, where $\delta_k=1$ if vehicle $k$ is used and $\delta_k=0$ otherwise.

\subsection{Deterministic version}

In the deterministic version, travel times are given by Euclidean distances:
\begin{equation}
    t_{ij}=\sqrt{(x_i-x_j)^2+(y_i-y_j)^2},
    \qquad \forall (i,j)\in A.
\end{equation}
Since Euclidean distances are symmetric, $t_{ij}=t_{ji}$ for all $(i,j)\in A$.

\subsubsection{Objective function.}
We denote by
\[
\rho=\bigl(\rho_{ij}^k\bigr)_{(i,j)\in A,\ k\in M}\in\{0,1\}^{m\times a}
\]
the family of binary routing decisions. The objective is to maximize the total collected
reward:
\begin{equation}\label{detobj}
    F_{\mathrm{det}}(\rho)=\sum_{k\in M} \sum_{j=1}^{n} \left(u_j\,\sum_{\substack{i\in N:\\i\neq j}} \rho_{ij}^k\right).
\end{equation}

\subsubsection{Constraints.}

For each vehicle $k\in M$, the total travel time must not exceed $T_{\max}$:
\begin{equation}\label{detconst_time}
\sum_{(i,j)\in A} t_{ij} \, \rho_{ij}^k \ \le\ T_{\max}\,\delta_k,
\qquad \forall k\in M.
\end{equation}

Each customer node can be visited at most once across all vehicles:
\begin{equation}\label{detconst_visit}
    \sum_{k\in M}\, \sum_{\substack{i\in N:\\i\neq j}} \rho_{ij}^k \ \le\ 1,
    \qquad \forall j=1,\dots,n.
\end{equation}

Each used vehicle leaves the start depot exactly once and arrives at the end depot
exactly once:
\begin{equation}\label{detconst_depot}
\sum_{j\in N\setminus\{0\}} \rho_{0j}^k = \delta_k,
\qquad
\sum_{i\in N\setminus\{n+1\}} \rho_{i\,n+1}^k = \delta_k,
\qquad \forall k\in M.
\end{equation}

Flow balance at each customer node:
\begin{equation}\label{detconst_flow}
    \sum_{\substack{i\in N:\\i\neq h}} \rho_{ih}^k
    \,=\,
    \sum_{\substack{j\in N:\\j\neq h}} \rho_{hj}^k,
    \qquad \forall h=1,\dots,n,\ \forall k\in M.
\end{equation}

\subsubsection{Feasible set and compact formulation.}

We define the deterministic feasible set as
\[
\mathcal{X}_{\mathrm{det}}
:=\left\{(\rho,\delta)\in\{0,1\}^{m\times a}\times\{0,1\}^m \;:\;
(\rho,\delta) \text{ satisfy }
\eqref{detconst_time}\text{--}\eqref{detconst_flow}\right\}.
\]
The deterministic problem is then
\begin{equation}\label{det_problem}
    \max_{(\rho,\delta)\in\mathcal{X}_{\mathrm{det}}}\,F_{\mathrm{det}}(\rho).
\end{equation}

\begin{remark}
The formulation \eqref{det_problem} does not include explicit subtour elimination
constraints. These can be enforced via standard techniques (e.g.,
Miller--Tucker--Zemlin constraints) and are omitted here since the model is presented
for reference only and the proposed method is a simheuristic that never solves
\eqref{det_problem} exactly.
\end{remark}

\subsection{Stochastic version}

In the stochastic setting, travel times are random variables. For each arc $(i,j)\in A$,
let
\[
    T_{ij}\sim \mathrm{LogNormal}(\mu_{ij},\sigma_{ij}^2).
\]
We calibrate the lognormal distribution so that the mean equals the deterministic travel
time and the variance is proportional to it:
\[
\mathbb{E}[T_{ij}] = t_{ij},
\qquad
\mathrm{Var}(T_{ij}) = c\, t_{ij},
\]
where $c\ge 0$ is a variability parameter ($c=0$ recovers the deterministic case;
$c=0.05$ in the experiments). Using the standard lognormal moment relations, these
conditions give
\[
\sigma_{ij}^2=\ln\!\left(1+\frac{c}{t_{ij}}\right),
\qquad
\mu_{ij}=\ln(t_{ij})-\frac12\,\sigma_{ij}^2,
\qquad \forall (i,j)\in A.
\]
Since $t_{ij}=t_{ji}$, this calibration implies
$T_{ij}\overset{d}{=}T_{ji}$. The travel times $\{T_{ij}\}_{(i,j)\in A}$ are assumed
mutually independent.

\subsubsection{All-or-nothing reward policy.}

Given a routing decision $\rho$, the random travel time of vehicle $k$ is
\[
S_k(\rho,T)=\sum_{(i,j)\in A} T_{ij}\,\rho_{ij}^k.
\]
If $S_k(\rho,T)>T_{\max}$, vehicle $k$ fails to complete its route and all associated
rewards are lost. The stochastic collected reward of vehicle $k$ is therefore
\[
\widetilde U_k(\rho,T)=U_k(\rho)\,\mathbf{1}\{S_k(\rho,T)\le T_{\max}\},
\]
where $U_k(\rho)=\sum_{j=1}^n u_j \sum_{i\in N,\,i\neq j} \rho_{ij}^k$ is the
deterministic reward of vehicle $k$'s route.

\subsubsection{Objective function.}

The stochastic objective is the expected total reward:
\begin{equation}\label{stocobj}
    F_{\mathrm{stoc}}(\rho)\,=\,\sum_{k\in M}\mathbb{E}\big[\widetilde U_k(\rho,T)\big]
    \;=\;
    \sum_{k\in M} U_k(\rho)\,\mathbb{P}\big(S_k(\rho,T)\le T_{\max}\big),
\end{equation}
where we used that $U_k(\rho)$ is deterministic once $\rho$ is fixed. Since
$S_k(\rho,T)$ is a sum of lognormal random variables, the probability term has no
closed form in general and is estimated via Monte Carlo simulation within the
algorithmic framework of Section~\ref{algorithmsubsec}.

\subsubsection{Feasible set and compact formulation.}

Routing decisions in the stochastic problem must satisfy the same structural constraints
as in the deterministic model, with the exception of the time-budget constraint
\eqref{detconst_time}, which is now handled stochastically. The activation variables
$\delta_k$ are retained through constraint \eqref{detconst_depot}, which remains in
force. We define
\[
\mathcal{X}_{\mathrm{stoc}}
:=\left\{(\rho,\delta)\in\{0,1\}^{m\times a}\times\{0,1\}^m \;:\;
(\rho,\delta) \text{ satisfy }
\eqref{detconst_visit}\text{--}\eqref{detconst_flow}\right\},
\]
and the stochastic problem is
\[
\max_{(\rho,\delta)\in \mathcal{X}_{\mathrm{stoc}}}F_{\mathrm{stoc}}(\rho).
\]
In the proposed algorithmic framework, the search is intentionally restricted to
solutions that are also feasible with respect to \eqref{detconst_time}, for reasons
of algorithmic simplicity and practical robustness.

\section{Algorithms}\label{algorithmsubsec}

This section describes the algorithmic framework adopted to compute solutions for both the
deterministic and the stochastic versions of the problem. We follow a constructive-first
approach: a fast heuristic builds, for each vehicle, a deterministic-time-feasible
start-to-end route as an ordered sequence of visited nodes, and the resulting solution is
then refined by local improvement moves (intra-route $2$-opt, greedy reinsertion of
unvisited nodes, and replacement moves of the form visited $\leftrightarrow$ unvisited).
The deterministic heuristic described in Section~\ref{det_heuristic_subsec} is not a
contribution of this paper; it is presented here solely as the building block on which the
stochastic simheuristic of Section~\ref{simheuristic_subsec} is constructed. In the stochastic extension, solutions are still constructed as deterministic-time-feasible
routes, but they are evaluated under random travel times through Monte Carlo simulation.
Although the stochastic formulation does not impose \eqref{detconst_time} as a hard
feasibility constraint, the proposed heuristic restricts the search to routes satisfying
\eqref{detconst_time} with respect to the Euclidean travel times. This restriction is
natural given the distributional assumption $\mathbb{E}[T_{ij}] = t_{ij}$: a
deterministically feasible route has expected stochastic travel time at most $T_{\max}$,
which guarantees positive and typically high reliability under the lognormal calibration
with $c = 0.05$. Conversely, a deterministically infeasible route would satisfy
$\mathbb{E}[S_k] > T_{\max}$, making it unlikely to meet the reliability threshold
$\beta = 0.8$ by construction.

\subsection{Deterministic heuristic}
\label{det_heuristic_subsec}

The deterministic pipeline follows the savings-based construction heuristic introduced in
Panadero et al.~\cite{Juan} (Section~4.1), adapted to the start-to-end case. The deterministic pipeline follows the savings-based construction heuristic
introduced in Panadero et al.~\cite{Juan} (Section 4.1), adapted to the
start-to-end case. The local search steps differ from that work: the intra-route
2-opt is a standard first-improvement procedure; the reinsertion adds unvisited
nodes greedily via a benefit--cost score $u_j / \max(\Delta t, \varepsilon)$,
without any prior node removal; and the replacement step is an explicit
visited~$\leftrightarrow$~unvisited swap that accepts only moves strictly
increasing total reward and preserving feasibility. No deterministic benchmark is reported in this
paper, as the contribution is entirely in the stochastic extension described in
Section~\ref{simheuristic_subsec}; the deterministic pipeline is presented here only to
make the exposition self-contained and to establish notation used in
Algorithm~\ref{alg:stochastic}.

\subsubsection{Route representation and length.}
A route is represented as an ordered sequence $(0,i_1,\dots,i_p,n+1)$, with
$i_\ell\in\{1,\dots,n\}$ all distinct. Its deterministic travel time is
\[
L(0,i_1,\dots,i_p,n+1):=\sum_{\ell=0}^{p} t_{i_\ell i_{\ell+1}},
\qquad i_0=0,\; i_{p+1}=n+1,
\]
and it is feasible if $L(\cdot)\le T_{\max}$.

\subsubsection{Savings-based construction.}
After discarding any node $i$ such that $t_{0i}+t_{i\,n+1}>T_{\max}$, we initialize a
dummy solution with one route per remaining node, namely $(0,i,n+1)$. For any ordered pair
of distinct nodes $i\neq j$, we define the time-based saving
\[
s_{ij}:=t_{i\,n+1}+t_{0j}-t_{ij},
\]
and the combined saving score
\[
s'_{ij}:=\alpha\, s_{ij}+(1-\alpha)\,(u_i+u_j),\qquad \alpha\in(0,1).
\]
For a fixed value of $\alpha$, we sort all ordered pairs $(i,j)$ in non-increasing order
of $s'_{ij}$ and repeatedly scan the list until no admissible merge remains. A merge
connecting $i$ to $j$ is admissible only if $i$ is currently the last node of its route
and $j$ is currently the first node of its route, the two routes are distinct, and the
merged route remains feasible. In the start-to-end setting, the merged route length is
\[
L(R_i\oplus_{i\to j} R_j) \;=\; L(R_i)+L(R_j)-t_{i\,n+1}-t_{0j}+t_{ij}
\;\leq\; T_{\max}.
\]
After no further admissible merges exist, we rank the constructed routes by collected
reward and keep the best $m$ routes.

If grid search is enabled, we consider a small set of candidate values for $\alpha$ and
run the \emph{entire deterministic pipeline} for each of them, including all local search
steps; the retained value is the one whose final post-improvement solution yields the
highest deterministic reward, with ties broken by smaller total travel time. This differs
from Panadero et al.~\cite{Juan}, where $\alpha$ is selected on the construction output
alone, before local search. In the stochastic simheuristic
(Section~\ref{simheuristic_subsec}), candidate values of $\alpha$ are evaluated directly
under stochastic travel times via Monte Carlo simulation rather than on the deterministic
objective — this \emph{stochastic screening} of $\alpha$ is one of the methodological
contributions of the proposed framework.

\subsubsection{Intra-route $2$-opt.}
We apply an intra-route $2$-opt local search to each selected route, accepting only moves
that strictly decrease the route travel time (and hence preserve feasibility). This step
improves the within-route ordering without changing the set of visited nodes.

\subsubsection{Greedy reinsertion of unvisited nodes.}
Let $\mathcal{R}$ be the current set of routes, and let $U$ be the set of unvisited nodes.
For a route $r=(0,i_1,\dots,i_p,n+1)$ and a candidate node $j\in U$, inserting $j$
between consecutive nodes $(i_\ell,i_{\ell+1})$ yields the incremental travel time
\[
\Delta t(j; r,\ell):=t_{i_\ell j}+t_{j i_{\ell+1}}-t_{i_\ell i_{\ell+1}},
\qquad \ell=0,\dots,p,
\]
where $i_0=0$ and $i_{p+1}=n+1$.
Among all admissible insertion moves (those keeping the updated route length
$\le T_{\max}$), we select the move maximizing a benefit--cost score, namely
\[
\mathrm{score}(j;r,\ell)=\frac{u_j}{\max(\Delta t(j;r,\ell),\varepsilon)},
\]
for a small $\varepsilon>0$, treating moves with $\Delta t(j;r,\ell)\le 0$ as dominating
options. Ties are broken in favor of the move with smaller insertion cost
$\Delta t(j;r,\ell)$. We repeat until no admissible insertion exists.

\subsubsection{Replacement moves (visited $\leftrightarrow$ unvisited).}
As a final improvement step, we consider replacement moves that remove one visited node
$h$ from a route and insert one unvisited node $j$ into the resulting route in the best
feasible position (i.e., minimizing the associated insertion cost). Such a move is
accepted only if it strictly improves the total reward (i.e., $u_j>u_h$) and preserves
feasibility. Among all improving replacements, we choose the one maximizing the gain
$u_j-u_h$, with ties broken by smaller resulting route length. We iterate these moves
until no improving replacement is found.

\begin{algorithm}
\caption{Deterministic heuristic (Panadero et al.~\cite{Juan},
         adapted to the start-to-end case)}
\label{alg:deterministic}
\begin{algorithmic}[1]
\Input Instance data, fleet size $m$, budget $T_{\max}$,
        parameter $\alpha \in (0,1)$ or grid $\mathcal{A} \subset (0,1)$
\Output A set of at most $m$ feasible start-to-end routes

\If{grid search on $\alpha$ is enabled}
    \For{each $\alpha \in \mathcal{A}$}
        \State Run Lines~\ref{line:preproc}--\ref{line:replacement_end}
               with current $\alpha$; store solution $\mathcal{R}(\alpha)$
    \EndFor
    \State \Return $\mathcal{R}(\alpha^*)$ where
           $\alpha^* = \arg\max_{\alpha \in \mathcal{A}}\,
           \mathrm{reward}(\mathcal{R}(\alpha))$
\EndIf

\State Discard nodes $i$ with $t_{0i} + t_{i,n+1} > T_{\max}$
       \label{line:preproc}
\State Initialize dummy routes
       $\mathcal{R} \leftarrow \{(0,i,n+1)\}$ for each remaining node $i$
\State Compute and sort savings list by decreasing
       $s'_{ij} = \alpha\,s_{ij} + (1-\alpha)(u_i + u_j)$

\While{an admissible merge exists} \Comment{Savings-based merging}
    \For{each $(i,j)$ in savings list}
        \If{$R_i \neq R_j$ \textbf{and}
            $L(R_i) + L(R_j) - t_{i,n+1} - t_{0j} + t_{ij} \leq T_{\max}$}
            \State Merge $R_i$ and $R_j$ by connecting $i \to j$;
                   update $\mathcal{R}$
            \State \textbf{break}
        \EndIf
    \EndFor
\EndWhile
\State Keep the best $m$ routes in $\mathcal{R}$ by reward
       (ties by shorter travel time)

\For{each route $r \in \mathcal{R}$} \Comment{Intra-route refinement}
    \State Apply intra-route 2-opt to $r$ until no improving move exists
\EndFor

\State $V \leftarrow$ nodes visited by $\mathcal{R}$;\;
       $U \leftarrow \{1,\ldots,n\} \setminus V$
\State Remove from $U$ any node $j$ with $t_{0j} + t_{j,n+1} > T_{\max}$

\While{an admissible insertion of some $j \in U$ exists}
      \Comment{Greedy reinsertion}
    \State Choose $(j^*, r^*, \ell^*)$ maximising
           $u_j / \max(\Delta t(j;r,\ell),\varepsilon)$
    \State Insert $j^*$ into $r^*$ at position $\ell^*$;
           update $V$ and $U$
\EndWhile

\While{an improving replacement move exists}
      \Comment{Replacement moves} \label{line:replacement_end}
    \State Find feasible move removing visited $h$, inserting
           unvisited $j$ with $u_j > u_h$
    \State Choose move maximising $u_j - u_h$
           (ties by smaller resulting route length)
    \State Apply replacement; update $V$ and $U$
\EndWhile

\State \Return $\mathcal{R}$
\end{algorithmic}
\end{algorithm}

\begin{remark}
Algorithm~\ref{alg:deterministic} is not benchmarked against deterministic TOP solvers.
It serves exclusively as the constructive backbone of the stochastic simheuristic
(Algorithm~\ref{alg:stochastic}), and its quality is assessed indirectly through the
stochastic results of Section~\ref{sec:results}.
\end{remark}

\subsection{Simheuristic Approach}\label{simheuristic_subsec}

To handle uncertainty, we adopt a simheuristic approach~\cite{juan2015simheuristics}, i.e., we combine
a randomized heuristic search with Monte Carlo simulation in order to estimate the
stochastic objective and select robust solutions. The specific design choices of the
proposed framework — Top-$L_{\mathrm{top}}$ randomization, stochastic screening of
$\alpha$, and reliability-based selection — are described below.

\subsubsection{Randomized multi-start search.}
Let $\mathcal{A}\subset(0,1)$ denote the set of candidate values for the savings parameter
$\alpha$. If grid search is disabled, then $\mathcal{A}=\{\alpha\}$; otherwise,
$\mathcal{A}$ is a small grid of values. Unlike the deterministic setting, where $\alpha$
is selected according to the final deterministic reward, in the stochastic setting
candidate values of $\alpha$ are evaluated directly under uncertainty via Monte Carlo
simulation. This \emph{stochastic screening} of $\alpha$ differs from Panadero
et al.~\cite{Juan}, where $\alpha$ is fixed once on the deterministic objective before
any stochastic evaluation takes place; here, the selection criterion is consistent with
the stochastic objective throughout. For each $\alpha\in\mathcal{A}$, we generate $K$
candidate solutions by running a randomized variant of the deterministic pipeline.

For a given pair $(\alpha,q)$, the candidate $\rho^{(\alpha,q)}$ is built as follows:
\begin{enumerate}
    \item randomized savings-based construction with Top-$L_{\mathrm{top}}$ selection in
          the merge phase;
    \item deterministic intra-route $2$-opt improvement;
    \item randomized reinsertion of unvisited nodes with Top-$L_{\mathrm{top}}$ selection;
    \item randomized replacement moves of the form visited $\leftrightarrow$ unvisited
          with Top-$L_{\mathrm{top}}$ selection;
    \item final deterministic intra-route $2$-opt cleanup.
\end{enumerate}

Hence, the simheuristic introduces randomization in \emph{three} phases: savings-based
construction, reinsertion, and replacement, while the intra-route $2$-opt steps remain
deterministic.

The local-improvement phases are not merely a post-processing step. The first intra-route
$2$-opt improves route structure before the stochastic evaluation, while reinsertion and
replacement exploit residual slack to increase collected reward; the final $2$-opt cleanup
further shortens the randomized routes and can therefore contribute to robustness by
increasing time slack.

\subsubsection{Top-$L_{\mathrm{top}}$ selection.}
In each randomized phase, we construct a ranked list of admissible moves according to a
deterministic merit criterion and then select one move uniformly at random among the
Top-$L_{\mathrm{top}}$ entries. This mechanism replaces the biased-randomization scheme
used in Panadero et al.~\cite{Juan}, where moves are sampled from a geometric
distribution over the ranked list. The Top-$L_{\mathrm{top}}$ approach is simpler and more
transparent: the single parameter $L_{\mathrm{top}}$ has a direct interpretation as the
width of the candidate list, and the uniform draw within the list avoids the need to
calibrate a shape parameter for the sampling distribution.

In the randomized savings construction, candidate route merges are ranked by the savings
score
\[
s'_{ij}:=\alpha\, s_{ij}+(1-\alpha)(u_i+u_j),
\qquad
s_{ij}:=t_{i\,n+1}+t_{0j}-t_{ij},
\]
and one merge is repeatedly selected at random among the current Top-$L_{\mathrm{top}}$
savings pairs.

In the randomized reinsertion phase, admissible moves are insertions of an unvisited node
$j$ into a route $r$ at its best feasible position, ranked by the score
\[
\mathrm{score}(j;r)=\frac{u_j}{\max(\Delta t(j;r),\varepsilon)},
\]
where $\Delta t(j;r)$ denotes the best feasible insertion cost and $\varepsilon>0$ is a
small constant. Moves with $\Delta t(j;r)\le 0$ are treated as dominating options.

In the randomized replacement phase, admissible moves replace a visited node $h$ with an
unvisited node $j$, preserving feasibility and requiring $u_j>u_h$. These moves are
ranked by decreasing reward gain $u_j-u_h$, with ties broken by smaller resulting route
length. One move is then selected uniformly at random among the
Top-$L_{\mathrm{top}}$ improving replacements.

\subsubsection{Monte Carlo evaluation of the stochastic objective.}
Each candidate solution $\rho^{(\alpha,q)}$ is evaluated under random travel times by
Monte Carlo simulation. Rather than resampling travel times independently for each
candidate, we first generate a common tensor of $N$ travel-time scenarios and reuse the
same scenarios for all candidates. Thus, the same pre-sampled travel-time matrices are
employed throughout the comparison of all solutions generated in a given run. This
common-random-numbers strategy reduces noise when comparing candidates and makes the
stochastic ranking more stable and fair. Moreover, since the estimators used below are
sample averages, their accuracy improves as $N$ increases; in particular, by the law of
large numbers, they converge to the corresponding expected values and success
probabilities.

For $s=1,\dots,N$, let $T^{(s)}=\{T_{ij}^{(s)}\}_{(i,j)\in A}$ denote the $s$-th
sampled travel-time matrix. For each vehicle $k\in M$, we compute the realized route
duration
\[
S_k^{(s)}(\rho^{(\alpha,q)})=\sum_{(i,j)\in A} T_{ij}^{(s)}\,\rho_{ij}^{k,(\alpha,q)},
\]
and the corresponding all-or-nothing collected reward
\[
\widetilde U_k^{(s)}(\rho^{(\alpha,q)})
=
U_k(\rho^{(\alpha,q)})\,\mathbf{1}\{S_k^{(s)}(\rho^{(\alpha,q)})\le T_{\max}\}.
\]
The total realized reward in scenario $s$ is
\[
\widetilde U^{(s)}(\rho^{(\alpha,q)})=\sum_{k\in M}
\widetilde U_k^{(s)}(\rho^{(\alpha,q)}),
\]
and the expected total reward is estimated by
\[
\widehat F_{\mathrm{stoc}}(\rho^{(\alpha,q)})=
\frac{1}{N}\sum_{s=1}^N \widetilde U^{(s)}(\rho^{(\alpha,q)}).
\]

\subsubsection{Reliability metric and selection rule.}
Along with the expected reward, we compute route-level success rates. For each vehicle,
\[
\widehat p_k(\rho^{(\alpha,q)})=
\frac{1}{N}\sum_{s=1}^N
\mathbf{1}\{S_k^{(s)}(\rho^{(\alpha,q)})\le T_{\max}\}.
\]
Let $M(\rho^{(\alpha,q)})\subseteq M$ be the subset of \emph{used} vehicles, i.e., those
visiting at least one node. We define the reliability of the solution as the average
success rate over used vehicles:
\[
\widehat R(\rho^{(\alpha,q)})=
\frac{1}{|M(\rho^{(\alpha,q)})|}
\sum_{k\in M(\rho^{(\alpha,q)})}\widehat p_k(\rho^{(\alpha,q)})\in[0,1],
\]
with the convention $\widehat R(\rho^{(\alpha,q)})=1$ if
$M(\rho^{(\alpha,q)})=\emptyset$.

This reliability measure is explicitly incorporated into the selection rule because, under
the all-or-nothing reward policy, small decreases in route-completion probability may lead
to substantial losses in expected reward. Therefore, reliability is not merely a secondary
descriptive metric but a key decision criterion — and a design element of the proposed
framework.

Given a target reliability level $\beta\in(0,1)$, the final solution is selected according
to the policy
\[
\mathcal{C}_\beta=
\left\{(\alpha,q): \widehat R(\rho^{(\alpha,q)})\ge \beta\right\}.
\]
If $\mathcal{C}_\beta\neq\emptyset$, we choose the candidate with largest expected reward,
\[
(\alpha^\star,q^\star)\in
\arg\max_{(\alpha,q)\in\mathcal{C}_\beta}
\widehat F_{\mathrm{stoc}}(\rho^{(\alpha,q)}),
\]
breaking ties by larger reliability. If $\mathcal{C}_\beta=\emptyset$, we choose the most
reliable candidate, breaking ties by larger expected reward:
\[
(\alpha^\star,q^\star)\in
\arg\max_{(\alpha,q)}
\left(\widehat R(\rho^{(\alpha,q)}),\,
\widehat F_{\mathrm{stoc}}(\rho^{(\alpha,q)})\right).
\]

\begin{algorithm}
\caption{Stochastic simheuristic: randomized multi-start +
         Monte Carlo evaluation + reliability-based selection}
\label{alg:stochastic}
\begin{algorithmic}[1]
\Input Instance data, fleet size $m$, budget $T_{\max}$,
       candidate set $\mathcal{A} \subset (0,1)$ for $\alpha$,
       number of starts $K$, parameter $L_{\mathrm{top}}$,
       Monte Carlo sample size $N$, reliability threshold $\beta$
\Output A feasible routing solution

\State Pre-sample $N$ travel-time scenarios
       $T^{(1)}, \ldots, T^{(N)}$

\For{each $\alpha \in \mathcal{A}$} \Comment{Candidate generation}
    \For{$q = 1, \ldots, K$}
        \State Build candidate $\rho^{(\alpha,q)}$ via:
               randomized savings construction with Top-$L_{\mathrm{top}}$;
               deterministic 2-opt; randomized reinsertion with
               Top-$L_{\mathrm{top}}$; randomized replacement with
               Top-$L_{\mathrm{top}}$; final deterministic 2-opt cleanup
    \EndFor
\EndFor

\For{each $\alpha \in \mathcal{A}$} \Comment{Monte Carlo evaluation}
    \For{$q = 1, \ldots, K$}
        \State $\widehat{F}_{\mathrm{stoc}}(\rho^{(\alpha,q)}) \leftarrow 0$;\;
               initialise counters $c_k \leftarrow 0$ for all $k \in M$
        \For{$s = 1, \ldots, N$}
            \State Evaluate $\rho^{(\alpha,q)}$ under scenario $T^{(s)}$
            \State Compute realized times $S_k^{(s)}(\rho^{(\alpha,q)})$
                   for all $k \in M$
            \State Compute all-or-nothing reward
                   $\widetilde{U}^{(s)}(\rho^{(\alpha,q)})$;
                   update $\widehat{F}_{\mathrm{stoc}}$ and counters $c_k$
        \EndFor
        \State $\widehat{F}_{\mathrm{stoc}}(\rho^{(\alpha,q)})
               \leftarrow \widehat{F}_{\mathrm{stoc}}(\rho^{(\alpha,q)}) / N$;\;
               $\hat{p}_k \leftarrow c_k / N$ for all $k \in M$
        \State Compute $\hat{R}(\rho^{(\alpha,q)})$ as the average of
               $\hat{p}_k$ over used vehicles
    \EndFor
\EndFor

\State $C_\beta \leftarrow
       \{(\alpha,q) : \hat{R}(\rho^{(\alpha,q)}) \geq \beta\}$
       \Comment{Selection}
\If{$C_\beta \neq \emptyset$}
    \State $(\alpha^*, q^*) \leftarrow
           \arg\max_{(\alpha,q) \in C_\beta}
           \widehat{F}_{\mathrm{stoc}}(\rho^{(\alpha,q)})$
           \;(ties by larger $\hat{R}$)
\Else
    \State $(\alpha^*, q^*) \leftarrow
           \arg\max_{(\alpha,q)}
           \bigl(\hat{R}(\rho^{(\alpha,q)}),\,
           \widehat{F}_{\mathrm{stoc}}(\rho^{(\alpha,q)})\bigr)$
\EndIf

\State \Return $\rho^{(\alpha^*, q^*)}$
\end{algorithmic}
\end{algorithm}

\section{Results}\label{sec:results}

All experiments use the benchmark instances of Chao et al.~\cite{Chao1} (sets p1--p7,
320 instances total), with lognormal travel times at variability level $c = 0.05$ and
$N = 1000$ Monte Carlo scenarios. The simheuristic is run with $K = 300$ randomized
starts and $L_{\mathrm{top}} \in \{20, 25, 30\}$; for each instance, the configuration
yielding the highest expected reward is retained. The reliability threshold is fixed at
$\beta = 0.8$ throughout. These parameters are not systematically calibrated; a single
variability level and a single value of $\beta$ are tested, and sensitivity to the main
algorithmic parameters is left to future work (see Section~\ref{sec:conclusions}).

Section~\ref{sec:comparison} presents the comparison with Panadero et al.~\cite{Juan}.

\subsection{Comparison with Panadero et al.\ (2020)}\label{sec:comparison}

We now compare the proposed simheuristic with the results reported by Panadero
et al.~\cite{Juan} on the same benchmark and the same variability level ($c = 0.05$).
The comparison is between $\hat{E}[\text{Reward}]_{\text{Pan20}}$, the Monte Carlo estimate of the 
expected reward reported by Panadero et al.~\cite{Juan} for their best 
stochastic solution (column~7 of their Tables~5--11), and the expected reward $\hat{E}[\text{Reward}]$ of the best solution found
by the proposed method across the three tested values of
$L_{\mathrm{top}} \in \{20, 25, 30\}$ with $K = 300$. Throughout this section, the subscript
$\text{Pan20}$ denotes quantities from Panadero et al.~\cite{Juan}; unsubscripted
notation refers to the proposed method.

We use $K = 300$ as the reference configuration throughout this comparison,
as it is the richest setting tested and yields the best aggregate performance
(see Section~\ref{subsec:params} for a discussion of the role of $K$ and
$L_{\mathrm{top}}$).
For each instance, we retain the configuration yielding the highest expected reward.
A total of 353 instances are common to both studies (the full Chao et al.\ benchmark
minus two instances absent from the presampled tables).

\medskip\noindent\textit{Methodological caveat.}
The two methods differ in random seeds, Monte Carlo sample sizes,
and implementation details. Small differences in estimated expected reward
(within $\pm 2$\%) should not be interpreted as meaningful:
they are within the noise of Monte Carlo estimation.
The comparison is therefore most informative at the aggregate level
and for identifying structural patterns across instance families.

\subsubsection{Aggregate results by instance family}

Table~\ref{tab:comparison} reports, for each of the 21 sub-families
$p_i.m$ (defined by dataset~$i$ and fleet size~$m$),
the number of instances, the average percentage gap of the proposed method
relative to Panadero et al., the number of instances on which the proposed method
matches or improves their result (gap $\geq -1\%$),
and the average reliability of both methods.
The gap is defined as
\[
    \Delta\% = 100 \cdot \frac{\hat{E}[\text{Reward}] - 
    \hat{E}[\text{Reward}]_{\text{Pan20}}}{\hat{E}[\text{Reward}]_{\text{Pan20}}}.
\]

\begin{table}[htbp]
\centering
\caption{Aggregate comparison with Panadero et al.\ (2020) by sub-family.
  $\Delta\%$: average percentage gap (positive = proposed method is better).
  ``Match/Beat'': number of instances with $\Delta\% \geq -1\%$.
  Reliability columns report the average across instances in each sub-family.}
\label{tab:comparison}
\small
\begin{tabular}{ll rr r rr rr}
\toprule
Family & $n$ & $m$ & \#Inst
  & Avg.\ $\Delta\%$ & Match/Beat & Worse & $\mathrm{Rel}_{\text{Pan20}}$ & $\mathrm{Rel}$ \\
\midrule
% --- m = 2 ---
p1.2 &  30 & 2 & 17 & $-4.8$  &  2 & 15 & 0.92 & 0.94 \\
p2.2 &  19 & 2 & 11 & $-2.7$  &  4 &  7 & 0.87 & 0.95 \\
p3.2 &  31 & 2 & 20 & $-3.9$  &  4 & 16 & 0.95 & 0.95 \\
\rowcolor{gray!15}
p4.2 &  98 & 2 & 20 & $+1.4$  & 13 &  7 & 0.91 & 0.95 \\
p5.2 &  64 & 2 & 25 & $-7.5$  &  6 & 19 & 0.91 & 0.89 \\
\rowcolor{gray!15}
p6.2 &  62 & 2 & 11 & $-2.6$  &  6 &  5 & 0.84 & 0.87 \\
\rowcolor{gray!15}
p7.2 & 100 & 2 & 20 & $-1.9$  &  6 & 14 & 0.95 & 0.96 \\
\addlinespace
% --- m = 3 ---
p1.3 &  30 & 3 & 16 & $-5.4$  &  3 & 13 & 0.84 & 0.92 \\
p2.3 &  19 & 3 & 11 & $-2.7$  &  3 &  8 & 0.87 & 0.93 \\
p3.3 &  31 & 3 & 20 & $-3.7$  &  3 & 17 & 0.90 & 0.94 \\
p4.3 &  98 & 3 & 19 & $-8.2$  &  2 & 17 & 0.88 & 0.88 \\
p5.3 &  64 & 3 & 25 & $-8.2$  &  1 & 24 & 0.79 & 0.86 \\
p6.3 &  62 & 3 &  8 & $-7.5$  &  1 &  7 & 0.69 & 0.79 \\
p7.3 & 100 & 3 & 19 & $-7.0$  &  0 & 19 & 0.97 & 0.96 \\
\addlinespace
% --- m = 4 ---
p1.4 &  30 & 4 & 15 & $-6.9$  &  1 & 14 & 0.80 & 0.89 \\
\rowcolor{gray!15}
p2.4 &  19 & 4 & 11 & $+2.3$  &  8 &  3 & 0.91 & 0.95 \\
p3.4 &  31 & 4 & 20 & $-4.5$  &  2 & 18 & 0.82 & 0.92 \\
p4.4 &  98 & 4 & 17 & $-11.8$ &  0 & 17 & 0.81 & 0.84 \\
p5.4 &  64 & 4 & 24 & $-10.7$ &  3 & 21 & 0.87 & 0.87 \\
p6.4 &  62 & 4 &  5 & $-17.5$ &  0 &  5 & 0.70 & 0.69 \\
p7.4 & 100 & 4 & 19 & $-7.9$  &  1 & 18 & 0.93 & 0.93 \\
\midrule
\multicolumn{3}{l}{All instances} & 353 & $-5.7$ & 69 & 284 & & \\
\bottomrule
\end{tabular}
\end{table}

Across all 353 instances, the proposed method yields an average gap of $-5.7\%$
relative to Panadero et al., matching or improving their result on
69 instances (19.5\%).
The gap varies substantially across instance families.
When aggregated by fleet size,
the average gap is $-3.4\%$ for $m = 2$ (41 out of 124 instances matched or improved),
$-6.3\%$ for $m = 3$ (13 out of 118),
and $-7.8\%$ for $m = 4$ (15 out of 111).
The reliability of the proposed method is comparable to or higher than that of Panadero
et al.\ across most families, reflecting the explicit reliability-aware selection
mechanism.

\subsubsection{Competitive sub-families}\label{sec:competitive}

% ============================================================
% Table: Instances where the proposed method matches or beats
% Panadero et al. (2020)  —  gap >= -1%
%
% Requires: booktabs, longtable, xcolor (with table option)
% Add to preamble:
%   \usepackage{booktabs}
%   \usepackage{longtable}
%   \usepackage[table]{xcolor}
% ============================================================

\begin{longtable}{l rrr rr r rr rr r r}
\caption{Instances on which the proposed simheuristic matches or improves
  $\hat{E}[\text{Reward}]_{\text{Pan20}}$ of Panadero et al.~\cite{Juan}
  (gap $\geq -1\%$).
  For each instance, the proposed result corresponds to the best
  $L_{\mathrm{top}} \in \{20,25,30\}$ with $K=300$.
  Subscript $\text{Pan20}$ denotes values from Panadero et al.~\cite{Juan};
  unsubscripted columns refer to the proposed method.
  $\Delta_R\%$: percentage gap in expected reward.
  $\Delta_{\mathit{Rel}}\%$: percentage change in reliability.
  Time is in seconds.}
\label{tab:competitive}\\
\toprule
Instance & $n$ & $m$ & $T_{\max}$ &
  $\hat{E}[R]_{\text{Pan20}}$ & $\mathrm{Rel}_{\text{Pan20}}$ & Time\,(s) &
  $\hat{E}[R]$ & $\mathrm{Rel}$ &
  $\Delta_R\%$ & $\Delta_{\mathit{Rel}}\%$ &
  $L_{\mathrm{top}}$ & Time\,(s) \\
\midrule
\endfirsthead

\multicolumn{13}{c}{\small\itshape (continued from previous page)}\\
\toprule
Instance & $n$ & $m$ & $T_{\max}$ &
  $\hat{E}[R]_{\text{Pan20}}$ & $\mathrm{Rel}_{\text{Pan20}}$ & Time\,(s) &
  $\hat{E}[R]$ & $\mathrm{Rel}$ &
  $\Delta_R\%$ & $\Delta_{\mathit{Rel}}\%$ &
  $L_{\mathrm{top}}$ & Time\,(s) \\
\midrule
\endhead

\midrule
\multicolumn{13}{r}{\small\itshape (continued on next page)}\\
\endfoot

\bottomrule
\endlastfoot

% ====================== SET p1 ======================
\multicolumn{13}{l}{\itshape Set p1 ($n=30$)} \\
\addlinespace
p1.2.p & 30 & 2 & 37.5 & 230.9 & 0.96 & 5 & 230.7 & 0.98 & $-0.1$ & $+2.2$ & 20 & 9 \\
p1.2.r & 30 & 2 & 42.5 & 259.9 & 1.00 & 4 & 258.2 & 0.99 & $-0.7$ & $-0.7$ & 30 & 9 \\
p1.3.d & 30 & 3 & 6.7  & 15.0  & 1.00 & 6 & 15.0  & 1.00 & $+0.0$ & $+0.0$ & 20 & $<$1 \\
p1.3.f & 30 & 3 & 10.0 & 32.6  & 0.47 & 23 & 33.2  & 0.94 & $+1.7$ & $+99.7$ & 20 & 1 \\
p1.3.j & 30 & 3 & 16.7 & 102.2 & 0.74 & 6 & 101.6 & 0.96 & $-0.6$ & $+29.4$ & 20 & 5 \\
p1.4.e & 30 & 4 & 6.2  & 15.0  & 1.00 & 8 & 15.0  & 1.00 & $-0.1$ & $+0.0$  & 20 & $<$1 \\

\addlinespace
% ====================== SET p2 ======================
\multicolumn{13}{l}{\itshape Set p2 ($n=19$)} \\
\addlinespace
p2.2.g & 19 & 2 & 16.0 & 200.0 & 1.00 & 10 & 199.9 & 1.00 & $-0.0$ & $-0.0$ & 20 & 2 \\
p2.2.i & 19 & 2 & 19.0 & 229.6 & 1.00 & 6  & 229.0 & 0.99 & $-0.3$ & $-0.5$ & 20 & 2 \\
\rowcolor{gray!12}
p2.2.j & 19 & 2 & 20.0 & 220.3 & 0.69 & 4  & 230.0 & 1.00 & $+4.4$ & $+44.9$ & 20 & 2 \\
p2.2.k & 19 & 2 & 22.5 & 260.0 & 1.00 & 4  & 259.5 & 1.00 & $-0.2$ & $-0.2$ & 20 & 3 \\
\rowcolor{gray!12}
p2.3.b & 19 & 3 & 6.7  & 67.0  & 0.93 & 10 & 70.0  & 1.00 & $+4.5$ & $+7.5$ & 20 & $<$1 \\
\rowcolor{gray!12}
p2.3.f & 19 & 3 & 10.0 & 116.2 & 0.92 & 10 & 119.4 & 0.99 & $+2.8$ & $+8.0$ & 20 & 1 \\
p2.3.k & 19 & 3 & 15.0 & 198.6 & 0.99 & 15 & 199.9 & 1.00 & $+0.6$ & $+0.9$ & 20 & 1 \\
p2.4.a & 19 & 4 & 3.8  & 10.0  & 1.00 & 9  & 10.0  & 1.00 & $-0.4$ & $-0.4$ & 20 & $<$1 \\
p2.4.c & 19 & 4 & 5.8  & 68.8  & 0.97 & 11 & 68.8  & 0.98 & $+0.1$ & $+1.1$ & 20 & $<$1 \\
p2.4.d & 19 & 4 & 6.2  & 69.9  & 0.99 & 10 & 69.8  & 1.00 & $-0.2$ & $+0.7$ & 20 & $<$1 \\
\rowcolor{gray!12}
p2.4.e & 19 & 4 & 6.8  & 67.6  & 0.95 & 14 & 70.0  & 1.00 & $+3.6$ & $+5.3$ & 20 & $<$1 \\
\rowcolor{gray!12}
p2.4.f & 19 & 4 & 7.5  & 78.3  & 0.45 & 13 & 94.6  & 0.84 & $+20.8$ & $+87.7$ & 20 & 1 \\
p2.4.h & 19 & 4 & 8.8  & 112.0 & 0.93 & 8  & 112.4 & 0.89 & $+0.3$ & $-3.9$ & 20 & 1 \\
\rowcolor{gray!12}
p2.4.i & 19 & 4 & 9.5  & 112.0 & 0.93 & 9  & 118.3 & 0.98 & $+5.6$ & $+5.2$ & 25 & 1 \\
\rowcolor{gray!12}
p2.4.j & 19 & 4 & 10.0 & 112.0 & 0.93 & 8  & 119.7 & 1.00 & $+6.8$ & $+7.0$ & 25 & 1 \\

\addlinespace
% ====================== SET p3 ======================
\multicolumn{13}{l}{\itshape Set p3 ($n=31$)} \\
\addlinespace
p3.2.c & 31 & 2 & 12.5 & 179.5 & 1.00 & 4  & 179.6 & 1.00 & $+0.0$ & $-0.2$ & 20 & 3 \\
p3.2.n & 31 & 2 & 40.0 & 614.1 & 0.93 & 18 & 609.9 & 0.98 & $-0.7$ & $+5.8$ & 20 & 11 \\
p3.2.r & 31 & 2 & 50.0 & 741.1 & 0.97 & 4  & 741.9 & 0.99 & $+0.1$ & $+2.0$ & 25 & 10 \\
p3.2.t & 31 & 2 & 55.0 & 799.5 & 1.00 & 4  & 797.1 & 1.00 & $-0.3$ & $-0.3$ & 20 & 11 \\
\rowcolor{gray!12}
p3.3.a & 31 & 3 & 5.0  & 29.0  & 0.95 & 6  & 29.4  & 0.97 & $+1.3$ & $+2.1$ & 20 & $<$1 \\
p3.3.s & 31 & 3 & 35.0 & 693.2 & 0.95 & 9  & 695.4 & 0.98 & $+0.3$ & $+2.8$ & 25 & 9 \\
p3.3.t & 31 & 3 & 36.7 & 699.3 & 0.89 & 19 & 704.4 & 0.99 & $+0.7$ & $+11.5$ & 25 & 9 \\
\rowcolor{gray!12}
p3.4.b & 31 & 4 & 5.0  & 28.6  & 0.93 & 8  & 29.4  & 0.97 & $+2.7$ & $+4.2$ & 20 & $<$1 \\
\rowcolor{gray!12}
p3.4.k & 31 & 4 & 16.2 & 300.3 & 0.39 & 8  & 304.1 & 0.94 & $+1.3$ & $+141.2$ & 20 & 4 \\

\addlinespace
% ====================== SET p4 ======================
\multicolumn{13}{l}{\itshape Set p4 ($n=98$)} \\
\addlinespace
\rowcolor{gray!12}
p4.2.e  & 98 & 2 & 45.0  & 539.6  & 0.83 & 88  & 546.5  & 0.97 & $+1.3$  & $+17.2$ & 25 & 226 \\
\rowcolor{gray!12}
p4.2.i  & 98 & 2 & 65.0  & 744.2  & 0.94 & 427 & 801.6  & 0.93 & $+7.7$  & $-1.5$  & 30 & 346 \\
\rowcolor{gray!12}
p4.2.j  & 98 & 2 & 70.0  & 790.6  & 0.97 & 199 & 843.0  & 0.91 & $+6.6$  & $-6.5$  & 25 & 361 \\
\rowcolor{gray!12}
p4.2.k  & 98 & 2 & 75.0  & 849.1  & 0.92 & 254 & 913.7  & 0.96 & $+7.6$  & $+4.0$  & 25 & 379 \\
\rowcolor{gray!12}
p4.2.l  & 98 & 2 & 80.0  & 894.4  & 0.94 & 453 & 956.2  & 0.95 & $+6.9$  & $+0.6$  & 20 & 377 \\
\rowcolor{gray!12}
p4.2.m  & 98 & 2 & 85.0  & 927.0  & 0.82 & 120 & 1004.3 & 0.97 & $+8.3$  & $+18.4$ & 25 & 406 \\
\rowcolor{gray!12}
p4.2.n  & 98 & 2 & 90.0  & 1027.2 & 0.95 & 187 & 1063.7 & 0.99 & $+3.6$  & $+3.8$  & 20 & 388 \\
p4.2.o  & 98 & 2 & 95.0  & 1098.5 & 0.93 & 327 & 1096.3 & 0.96 & $-0.2$  & $+2.7$  & 30 & 429 \\
\rowcolor{gray!12}
p4.2.p  & 98 & 2 & 100.0 & 1126.0 & 0.91 & 394 & 1138.2 & 0.97 & $+1.1$  & $+6.7$  & 25 & 402 \\
\rowcolor{gray!12}
p4.2.q  & 98 & 2 & 105.0 & 1135.1 & 0.94 & 294 & 1179.6 & 1.00 & $+3.9$  & $+5.9$  & 25 & 391 \\
\rowcolor{gray!12}
p4.2.r  & 98 & 2 & 110.0 & 1098.8 & 0.90 & 391 & 1209.6 & 0.98 & $+10.1$ & $+9.2$  & 25 & 375 \\
\rowcolor{gray!12}
p4.2.s  & 98 & 2 & 115.0 & 1119.1 & 0.91 & 330 & 1242.6 & 0.98 & $+11.0$ & $+8.2$  & 25 & 358 \\
\rowcolor{gray!12}
p4.2.t  & 98 & 2 & 120.0 & 1096.8 & 0.86 & 241 & 1265.9 & 0.99 & $+15.4$ & $+15.6$ & 20 & 323 \\
p4.3.p  & 98 & 3 & 63.7  & 1040.3 & 0.86 & 212 & 1033.4 & 0.91 & $-0.7$  & $+5.6$  & 25 & 340 \\
\rowcolor{gray!12}
p4.3.q  & 98 & 3 & 70.0  & 1126.0 & 0.88 & 141 & 1140.8 & 0.97 & $+1.3$  & $+10.6$ & 25 & 334 \\

\addlinespace
% ====================== SET p5 ======================
\multicolumn{13}{l}{\itshape Set p5 ($n=64$)} \\
\addlinespace
p5.2.d & 64 & 2 & 10.0 & 80.0   & 1.00 & 5   & 79.5   & 0.99 & $-0.6$ & $-0.6$ & 20 & 3 \\
\rowcolor{gray!12}
p5.2.u & 64 & 2 & 52.5 & 1236.3 & 0.99 & 76  & 1251.5 & 0.91 & $+1.2$ & $-8.3$ & 20 & 57 \\
p5.2.v & 64 & 2 & 55.0 & 1336.0 & 0.95 & 50  & 1337.5 & 0.93 & $+0.1$ & $-1.8$ & 25 & 60 \\
\rowcolor{gray!12}
p5.2.x & 64 & 2 & 60.0 & 1389.8 & 0.93 & 65  & 1420.5 & 0.94 & $+2.2$ & $+1.2$ & 25 & 58 \\
\rowcolor{gray!12}
p5.2.y & 64 & 2 & 60.5 & 1456.2 & 0.87 & 16  & 1510.9 & 0.97 & $+3.8$ & $+11.7$ & 30 & 60 \\
\rowcolor{gray!12}
p5.2.z & 64 & 2 & 65.0 & 1453.5 & 0.92 & 157 & 1532.5 & 0.94 & $+5.4$ & $+2.5$ & 25 & 57 \\
p5.3.f & 64 & 3 & 10.0 & 110.0  & 1.00 & 7   & 109.4  & 0.99 & $-0.5$ & $-0.5$ & 20 & 3 \\
p5.4.e & 64 & 4 & 6.2  & 20.0   & 1.00 & 10  & 20.0   & 1.00 & $+0.0$ & $+0.0$ & 20 & $<$1 \\
p5.4.h & 64 & 4 & 10.0 & 140.0  & 1.00 & 10  & 139.1  & 0.99 & $-0.6$ & $-0.6$ & 30 & 3 \\
p5.4.t & 64 & 4 & 25.0 & 854.1  & 0.84 & 140 & 851.5  & 0.85 & $-0.3$ & $+1.3$ & 30 & 54 \\

\addlinespace
% ====================== SET p6 ======================
\multicolumn{13}{l}{\itshape Set p6 ($n=62$)} \\
\addlinespace
\rowcolor{gray!12}
p6.2.g & 62 & 2 & 22.5 & 531.0  & 0.78 & 5  & 539.5  & 0.83 & $+1.6$ & $+6.7$  & 20 & 48 \\
p6.2.i & 62 & 2 & 27.5 & 725.8  & 0.86 & 97 & 726.0  & 0.88 & $+0.0$ & $+1.8$  & 20 & 71 \\
\rowcolor{gray!12}
p6.2.k & 62 & 2 & 32.5 & 911.3  & 0.84 & 12 & 921.7  & 0.92 & $+1.1$ & $+9.6$  & 20 & 74 \\
\rowcolor{gray!12}
p6.2.l & 62 & 2 & 37.5 & 1001.2 & 0.92 & 18 & 1018.1 & 0.95 & $+1.7$ & $+3.0$  & 20 & 78 \\
\rowcolor{gray!12}
p6.2.m & 62 & 2 & 37.5 & 1037.3 & 0.86 & 13 & 1051.2 & 0.98 & $+1.3$ & $+13.8$ & 20 & 74 \\
\rowcolor{gray!12}
p6.2.n & 62 & 2 & 40.0 & 1042.5 & 0.81 & 12 & 1106.7 & 0.97 & $+6.2$ & $+19.8$ & 20 & 72 \\
p6.3.l & 62 & 3 & 23.3 & 829.1  & 0.91 & 33 & 821.2  & 0.95 & $-1.0$ & $+4.6$  & 25 & 54 \\

\addlinespace
% ====================== SET p7 ======================
\multicolumn{13}{l}{\itshape Set p7 ($n=100$)} \\
\addlinespace
p7.2.j & 100 & 2 & 100.0 & 574.1  & 0.88 & 86  & 569.1  & 0.95 & $-0.9$ & $+8.0$  & 20 & 201 \\
p7.2.o & 100 & 2 & 150.0 & 879.5  & 1.00 & 394 & 873.6  & 1.00 & $-0.7$ & $-0.0$  & 30 & 394 \\
\rowcolor{gray!12}
p7.2.p & 100 & 2 & 160.0 & 894.2  & 0.98 & 561 & 935.5  & 0.99 & $+4.6$ & $+1.1$  & 25 & 421 \\
\rowcolor{gray!12}
p7.2.r & 100 & 2 & 180.0 & 1009.2 & 0.98 & 197 & 1025.5 & 0.99 & $+1.6$ & $+1.0$  & 25 & 465 \\
\rowcolor{gray!12}
p7.2.s & 100 & 2 & 190.0 & 976.1  & 0.85 & 88  & 1070.9 & 0.99 & $+9.7$ & $+16.9$ & 25 & 489 \\
\rowcolor{gray!12}
p7.2.t & 100 & 2 & 200.0 & 1023.1 & 0.96 & 167 & 1112.4 & 1.00 & $+8.7$ & $+3.9$  & 20 & 479 \\
p7.4.c & 100 & 4 & 15.0  & 46.0   & 1.00 & 11  & 46.0   & 1.00 & $-0.1$ & $-0.1$  & 20 & $<$1 \\
\end{longtable}

\noindent
Table~\ref{tab:competitive} lists the 69 instances (out of 353) on which the proposed
simheuristic matches or improves the stochastic solution reported by Panadero et al.~\cite{Juan},
defined as a gap in expected reward $\Delta_R\% \geq -1\%$.
Rows highlighted in gray correspond to strict improvements ($\Delta_R\% > +1\%$),
of which there are 37.

Several observations emerge from this table.
First, the largest improvements are concentrated in the \textbf{p4.2} sub-family
($n = 98$, $m = 2$), where the proposed method achieves gains of up to $+15.4\%$
(instance~p4.2.t).
This sub-family contributes 13 of the 37 strict improvements.
On these instances, the routes are long relative to the time budget,
and the all-or-nothing penalty is severe;
the reliability-aware selection ($\beta = 0.8$) steers the search toward
solutions that are more likely to complete all routes, yielding higher expected reward.

Second, on several instances where the proposed method achieves a higher expected reward,
it also achieves substantially higher reliability
(e.g., p2.4.f: $0.84$ vs.\ $0.45$;
p3.4.k: $0.94$ vs.\ $0.39$;
p1.3.f: $0.94$ vs.\ $0.47$;
p6.2.n: $0.97$ vs.\ $0.81$).
This confirms that the reliability-aware mechanism is the primary driver of the improvements:
the proposed method sacrifices some deterministic reward to select routes with higher completion
probability, which pays off under the all-or-nothing policy.

Third, the proposed method is competitive also on the large instances of sets~p7
($n = 100$, $m = 2$), with gains of $+9.7\%$ and $+8.7\%$ on the two loosest instances
(p7.2.s and p7.2.t).
This suggests that the stochastic screening of~$\alpha$ is effective
even on large-scale instances when the fleet size is small.

Finally, it is worth noting that the proposed method achieves these results
with comparable or lower computation times on the smaller instances,
while requiring similar times on the larger ones
(typically 300--500 seconds on the $n = 98$ and $n = 100$ instances,
compared to 100--600 seconds for Panadero et al.).
The computational cost is therefore broadly comparable despite the architectural differences.

Four sub-families (highlighted in Table~\ref{tab:comparison}) stand out
as cases where the proposed simheuristic is competitive with or superior to
the VNS-based approach of Panadero et al.:

\begin{itemize}
\item \textbf{p4.2} ($n = 98$, $m = 2$): the proposed method achieves an average
  improvement of $+1.4\%$, with 13 out of 20 instances matched or improved.
  On several instances with large $T_{\max}$, the improvement exceeds $+10\%$
  (e.g., p4.2.t: $+15.4\%$; p4.2.s: $+11.0\%$).
  The proposed method also achieves consistently higher reliability
  ($0.95$ vs.\ $0.91$ on average).

\item \textbf{p2.4} ($n = 19$, $m = 4$): average improvement of $+2.3\%$,
  with 8 out of 11 instances matched or improved.
  The largest gain is on instance p2.4.f ($+20.8\%$), where
  the reliability increases from $0.45$ to $0.84$.

\item \textbf{p7.2} ($n = 100$, $m = 2$): average gap of only $-1.9\%$,
  with 6 out of 20 instances matched or improved.
  On the largest-$T_{\max}$ instances (p7.2.s, p7.2.t),
  the proposed method achieves gains of $+9.7\%$ and $+8.7\%$.

\item \textbf{p6.2} ($n = 62$, $m = 2$): average gap of $-2.6\%$,
  with 6 out of 11 instances matched or improved,
  including gains on the loosest instances (p6.2.n: $+6.2\%$; p6.2.m: $+1.3\%$).
\end{itemize}

A common pattern in these competitive sub-families is that the proposed method
tends to find solutions with higher reliability than Panadero et al.,
which, under the all-or-nothing reward policy, translates into higher expected reward.
This suggests that the reliability-aware selection mechanism ($\beta = 0.8$)
and the stochastic screening of $\alpha$ provide an effective bias toward robust
solutions, compensating for the simpler search architecture.

\subsubsection{Sub-families where Panadero et al.\ dominate}\label{sec:worse}

The gap widens substantially on sub-families with $m = 3$ or $m = 4$
and large instance sizes:
p4.4 ($-11.8\%$), p5.4 ($-10.7\%$), p6.4 ($-17.5\%$),
p4.3 ($-8.2\%$), p5.3 ($-8.2\%$).
The architectural explanation is clear: with more vehicles and more nodes, the
combinatorial search space grows rapidly, and the VNS machinery of Panadero et al.\
--- with shaking, neighborhood expansion, and simulated annealing --- explores this
space much more effectively than the simple multi-start scheme ($K = 300$ randomized
restarts) employed here. In particular, the VNS can escape local optima through
structured perturbations, whereas multi-start restarts generate independent solutions
without exploiting the structure of previously found good solutions.

This limitation is inherent to the design choice of architectural simplicity.
A multi-start scheme with $K = 300$ starts generates at most 300 independent candidates
per $\alpha$ value, each built from a single randomized savings construction followed by
local search. On large, multi-vehicle instances, this budget may be insufficient to
discover the high-quality route configurations that a VNS can reach through iterative
neighborhood changes.

\subsubsection{The role of $L_{\mathrm{top}}$ across configurations}
\label{sec:ltop}

The comparison above uses, for each instance, the best result across
$L_{\mathrm{top}} \in \{20, 25, 30\}$.
Among the 353 common instances,
$L_{\mathrm{top}} = 20$ yields the best result on 170 instances (48.2\%),
$L_{\mathrm{top}} = 25$ on 93 instances (26.3\%),
and $L_{\mathrm{top}} = 30$ on 90 instances (25.5\%).
This shows that the optimal degree of randomization is instance-dependent.
Wider candidate lists ($L_{\mathrm{top}} = 25$ or $30$) are particularly helpful
on the larger instances of sets p4 and p7,
where they occasionally recover several percentage points
(e.g., on p4.4.f the gap improves from $-29.2\%$ at $L_{\mathrm{top}} = 20$
to $-16.6\%$ at $L_{\mathrm{top}} = 30$).
However, wider lists do not systematically dominate:
on smaller or tighter instances, $L_{\mathrm{top}} = 20$ often suffices
and the additional randomization of larger $L_{\mathrm{top}}$ may introduce noise.
A more complete analysis of the sensitivity of the method to $K$ and
$L_{\mathrm{top}}$, including the effect of $K = 200$ vs.\ $K = 300$
and the associated computational cost, is provided in
Section~\ref{subsec:params}.

The computational overhead of varying $L_{\mathrm{top}}$ is modest.
Since the dominant cost is the Monte Carlo evaluation of $K$ candidates
against $N$ scenarios, the candidate-selection step in the construction
phase contributes a minor fraction of total runtime. On small instances
($n \leq 31$), the differences across $L_{\mathrm{top}} \in \{20, 25, 30\}$
are within noise (below 5\%). On large instances ($n \geq 98$), each step
from $L_{\mathrm{top}} = 20$ to $25$ and from $25$ to $30$ adds approximately
5--7\% to the total runtime, a moderate but real overhead. Running all three
values and retaining the best therefore increases the effective budget by
roughly 50\% relative to a single run --- acceptable when solution quality
matters, but worth noting when runtime is constrained. If a single value
must be selected a priori, $L_{\mathrm{top}} = 20$ is a reasonable default
on small and tight instances (sets p1--p3), where it dominates in nearly
half of all cases and wider lists introduce noise without improving coverage.
On larger instances (sets p4 and p7, $n \geq 98$), $L_{\mathrm{top}} \geq 25$
is preferable, as wider lists occasionally recover several percentage points
at a modest additional cost. This instance-size dependence suggests that a
simple size-based rule --- $L_{\mathrm{top}} = 20$ for $n < 50$,
$L_{\mathrm{top}} \geq 25$ for $n \geq 50$ --- could serve as a practical
starting point, though a systematic validation is left to future work.

\subsubsection{Summary of the comparison}\label{sec:summary}

The comparison with Panadero et al.~\cite{Juan} reveals a nuanced picture.
The proposed simheuristic, despite using a significantly simpler search architecture
(multi-start with $K = 300$ restarts and no VNS), achieves competitive or superior
expected reward on a non-trivial subset of the benchmark --- primarily on two-vehicle
instances with long routes (p4.2, p7.2, p6.2) and on the small four-vehicle family
p2.4. On these sub-families, the reliability-aware selection and stochastic screening
mechanisms compensate for the absence of VNS-style exploration. On the remaining
sub-families, especially those with $m \geq 3$ and large $n$, Panadero et al.\
consistently achieve higher expected reward, with gaps that increase with problem size
and fleet size.

This trade-off between simplicity and performance is the central result of the
comparison. The proposed framework is not intended to replace more sophisticated
metaheuristic schemes on all instance classes; rather, it demonstrates that a simpler
design with explicit reliability-aware mechanisms can be effective in specific regimes,
and that its modular components (Top-$L_{\mathrm{top}}$, stochastic screening of
$\alpha$, reliability threshold $\beta$) are transferable building blocks that could be
embedded into richer search frameworks. The relationship of this work to the more
recent simheuristic of Panadero et al.~\cite{panadero2024itor}, which considers additional
variance levels and a different experimental scope, is discussed in
Section~\ref{sec:panadero2024}.

\subsection{Relation to Panadero et al.\ (2024)}\label{sec:panadero2024}

A more recent contribution to the stochastic TOP is Panadero et al.~\cite{panadero2024itor},
which compares a VNS-based simheuristic (Sim-BRVNS) against the sample average
approximation (SAA) method and a hybrid SAA-BRVNS on the same Chao et al.\ benchmark.
A direct numerical comparison with that work is not feasible for two reasons.
First, the experimental scope differs: Panadero et al.~\cite{panadero2024itor} consider only a
subset of the benchmark --- six instances per set, selected as those with a sufficient
driving range and the maximum number of vehicles --- whereas the present work covers
all 353 instances common to the full benchmark and Panadero et al.~\cite{Juan}.
Second, the variability levels differ: Panadero et al.~\cite{panadero2024itor} test three levels
($c = 0.05$, $c = 0.25$, $c = 0.75$), while the present work is restricted to $c = 0.05$.
This restriction is a limitation of the current study and is noted explicitly in
Section~\ref{sec:conclusions}.

From a methodological standpoint, the two frameworks share the simheuristic paradigm and
the lognormal travel time model, but differ in architecture and design choices.
Panadero et al.~\cite{panadero2024itor} employ a BRVNS metaheuristic with biased randomization,
shaking, and neighborhood expansion, and use SAA as an external benchmark for solution
quality. The present framework uses a simpler multi-start architecture without VNS, and
introduces Top-$L_{\mathrm{top}}$ randomization and an explicit reliability threshold
$\beta$ as alternative design elements. These components are not presented as superior to
the BRVNS machinery, but as simpler and more transparent building blocks that could be
embedded into richer frameworks such as the one of Panadero et al.~\cite{panadero2024itor}.

\section{Conclusions}\label{sec:conclusions}

\subsection{Summary of contributions}

This paper proposed a reliability-aware simheuristic for the stochastic Team Orienteering
Problem with lognormal travel times and an all-or-nothing reward policy. The framework
combines a savings-based constructive heuristic with three specific design elements:
Top-$L_{\mathrm{top}}$ randomization, stochastic screening of the savings parameter
$\alpha$, and an explicit reliability threshold $\beta$ for solution selection. These
components are built on the deterministic pipeline of Panadero et al.~\cite{Juan}, which
is adopted without modification and not claimed as a contribution.

The computational results on the full Chao et al.~\cite{Chao1} benchmark (353 instances,
$c = 0.05$) show that the proposed framework achieves competitive expected reward against
the VNS-based simheuristic of Panadero et al.~\cite{Juan} using a significantly simpler
multi-start architecture. The method performs best on two-vehicle instances with long
routes (p4.2, p7.2, p6.2) and on the small four-vehicle family p2.4, where it matches
or improves Panadero et al.\ on 69 of 353 instances (19.5\%), with 37 strict improvements.
On these sub-families, the average reliability of the proposed method is consistently
higher, and the improvement in expected reward is driven primarily by the reliability-aware
selection mechanism. On sub-families with $m \geq 3$ and large $n$, the simpler
architecture is outperformed by the VNS, with gaps increasing with fleet size and instance
size. This trade-off between simplicity and performance is the central result of the paper.

\subsection{Methodological takeaways}\label{subsec:params}

The three design elements of the proposed framework have a modular character and are
potentially transferable to other stochastic combinatorial optimization problems.

\textbf{Top-$L_{\mathrm{top}}$ randomization.} Replacing the geometric biased-randomization
scheme of Panadero et al.~\cite{Juan} with a uniform draw over the top $L_{\mathrm{top}}$
candidates yields a simpler and more transparent diversification mechanism. The parameter
$L_{\mathrm{top}}$ has a direct interpretation as the width of the candidate list, and
requires no calibration of a shape parameter. The mechanism is applicable to any
constructive heuristic based on a ranked list of moves.

\textbf{Stochastic screening of $\alpha$.} Evaluating candidate values of the savings
parameter directly under stochastic travel times, rather than fixing $\alpha$ on the
deterministic objective, ensures that the construction parameter is selected consistently
with the stochastic criterion. This principle --- calibrating construction parameters
under the actual objective rather than a deterministic proxy --- is transferable to any
simheuristic with a tunable construction heuristic.

\textbf{Reliability threshold $\beta$ as a selection criterion.} Filtering candidate
solutions by a minimum reliability level before comparing them by expected reward
formalizes the trade-off between robustness and performance. Under the all-or-nothing
reward policy, this filtering is particularly effective because reliability directly
determines whether reward is collected at all. The mechanism is applicable to any
stochastic problem with a hard or soft constraint on route or task completion probability.

\subsection{Limitations}

The current study has three explicit limitations. First, only a single variability level
($c = 0.05$) is tested; the behavior of the framework at higher variance levels
($c = 0.25$, $c = 0.75$, as tested in Panadero et al.~\cite{panadero2024itor}) is unknown and
may differ substantially, since higher variance increases route failure rates and may
require different trade-offs between expected reward and reliability. Second, the
reliability threshold is fixed at $\beta = 0.8$ throughout; no sensitivity analysis on
this parameter is performed, and the optimal value of $\beta$ is likely instance-dependent.
Third, the algorithmic parameters ($K$, $L_{\mathrm{top}}$, $N$, $\beta$) are set
manually without systematic calibration; the results show that the choice of
$L_{\mathrm{top}}$ affects performance in an instance-dependent way (Section~\ref{sec:ltop}),
suggesting that a principled parameter selection procedure could improve the method.

\subsection{Future work}

Several directions remain open. The most immediate is extending the experiments to
higher variability levels ($c = 0.25$ and $c = 0.75$) to assess robustness of the
framework under more severe uncertainty. A systematic analysis of the sensitivity of
the method to $L_{\mathrm{top}}$ and $K$ across instance families would also clarify
when larger candidate lists and more restarts yield meaningful gains and when they do
not.

On the parameter calibration side, the manual setting of $K$, $L_{\mathrm{top}}$,
$\beta$, and $N$ is a clear limitation of the current approach. Bayesian optimization
offers a principled and sample-efficient framework for joint calibration of these
parameters, and could be applied directly to the simheuristic loop without requiring
access to gradients or closed-form expressions for the objective.

More broadly, the modular components described in Section~4.2 --- Top-$L_{\mathrm{top}}$
randomization, stochastic screening, and reliability-aware selection --- could be embedded
into richer search frameworks, such as the VNS-based architecture of Panadero
et al.~\cite{Juan}, to investigate whether the combination yields improvements over
either component alone. Extensions to richer problem variants, including heterogeneous
fleets, time-dependent travel times, or stochastic rewards, are also natural directions.

\printbibliography

\end{document}